\documentclass[12pt]{amsart}
\usepackage{amssymb,latexsym}
\usepackage{array}
\usepackage{tikz-cd}
\usepackage{amsmath, amssymb, amsfonts, amstext, amsthm, amscd}
\usepackage[mathscr]{euscript}
\usepackage{enumerate}
\usepackage{tikz,color,soul}
\usepackage[Symbolsmallscale]{upgreek}
\usetikzlibrary{arrows.meta, positioning,arrows}
\usetikzlibrary{matrix,arrows}

\usepackage[english]{babel}
\usepackage{chngcntr}
\usepackage{textgreek}
\begin{document}

\thispagestyle{empty} 

\title[Equivariant $K$-theory of cellular toric bundles]{Equivariant $K$-theory of cellular toric bundles and related
  spaces}

\author{V. Uma}

\address{Department of Mathematics, Indian Institute of Technology Madras, Chennai-600036}
\email{vuma@iitm.ac.in}




\date{}

\begin{abstract}
  In this article we describe the equivariant and ordinary topological
  $K$-ring of a toric bundle with fiber a $T$-{\it cellular} toric
  variety. This generalizes the results in \cite{su} on $K$-theory of
  smooth projective toric bundles. We apply our results to describe
  the equivariant topological $K$-ring of a toroidal horospherical
  embedding.\end{abstract}

\maketitle

\def\theequation
  {\arabic{section}.\arabic{equation}}

\newcommand{\codim}{\mbox{{\rm codim}$\,$}}
\newcommand{\stab}{\mbox{{\rm stab}$\,$}}
\newcommand{\lr}{\mbox{$\rightarrow$}}

\newcommand{\be}{\begin{equation}}
\newcommand{\ee}{\end{equation}}

\newtheorem{guess}{Theorem}[section]
\newcommand{\bth}{\begin{guess}$\!\!\!${\bf }~}
\newcommand{\eeth}{\end{guess}}
\renewcommand{\bar}{\overline}
\newtheorem{propo}[guess]{Proposition}
\newcommand{\bpropo}{\begin{propo}$\!\!\!${\bf }~}
\newcommand{\epropo}{\end{propo}}

\newtheorem{assum}[guess]{Assumption}
\newcommand{\bas}{\begin{assum}$\!\!\!${\bf }~}
\newcommand{\eas}{\end{assum}}

\newtheorem{lema}[guess]{Lemma}
\newcommand{\blem}{\begin{lema}$\!\!\!${\bf }~}
\newcommand{\elem}{\end{lema}}

\newtheorem{defe}[guess]{Definition}
\newcommand{\bdefe}{\begin{defe}$\!\!\!${\bf }~}
\newcommand{\edefe}{\end{defe}}

\newtheorem{coro}[guess]{Corollary}
\newcommand{\bcor}{\begin{coro}$\!\!\!${\bf }~}
\newcommand{\ecor}{\end{coro}}

\newtheorem{rema}[guess]{Remark}
\newcommand{\brem}{\begin{rema}$\!\!\!${\bf }~\rm}
\newcommand{\erem}{\end{rema}}

\newtheorem{exam}[guess]{Example}
\newcommand{\beg}{\begin{exam}$\!\!\!${\bf }~\rm}
\newcommand{\eeg}{\end{exam}}

\newtheorem{notn}[guess]{Notation}
\newcommand{\bnot}{\begin{notn}$\!\!\!${\bf }~\rm}
\newcommand{\enot}{\end{notn}}

\newcommand{\ch}{{\mathcal H}}
\newcommand{\cf}{{\mathcal F}}
\newcommand{\cd}{{\mathcal D}}
\newcommand{\cR}{{\mathcal R}}
\newcommand{\cv}{{\mathcal V}}
\newcommand{\cn}{{\mathcal N}}
\newcommand{\lra}{\rightarrow}
\newcommand{\ra}{\rightarrow}
\newcommand{\blr}{\Big \rightarrow}
\newcommand{\da}{\Big \downarrow}
\newcommand{\ua}{\Big \uparrow}
\newcommand{\hra}{\mbox{{$\hookrightarrow$}}}
\newcommand{\rt}{\mbox{\Large{$\rightarrowtail$}}}
\newcommand{\dua}{\begin{array}[t]{c}
\Big\uparrow \\ [-4mm]
\scriptscriptstyle \wedge \end{array}}
\newcommand{\ctext}[1]{\makebox(0,0){#1}}
\setlength{\unitlength}{0.1mm}
\newcommand{\cl}{{\mathcal L}}
\newcommand{\cp}{{\mathcal P}}
\newcommand{\ci}{{\mathcal I}}
\newcommand{\bz}{\mathbb{Z}}
\newcommand{\cs}{{\mathcal s}}
\newcommand{\ce}{{\mathcal E}}
\newcommand{\ck}{{\mathcal K}}
\newcommand{\cz}{{\mathcal Z}}
\newcommand{\cg}{{\mathcal G}}
\newcommand{\cj}{{\mathcal J}}
\newcommand{\cc}{{\mathcal C}}
\newcommand{\ca}{{\mathcal A}}
\newcommand{\cb}{{\mathcal B}}
\newcommand{\cx}{{\mathcal X}}
\newcommand{\co}{{\mathcal O}}
\newcommand{\bq}{\mathbb{Q}}
\newcommand{\bt}{\mathbb{T}}
\newcommand{\bh}{\mathbb{H}}
\newcommand{\br}{\mathbb{R}}
\newcommand{\bl}{\mathbf{L}}
\newcommand{\wt}{\widetilde}
\newcommand{\im}{{\rm Im}\,}
\newcommand{\bc}{\mathbb{C}}
\newcommand{\bp}{\mathbb{P}}
\newcommand{\ba}{\mathbb{A}}
\newcommand{\spin}{{\rm Spin}\,}
\newcommand{\ds}{\displaystyle}
\newcommand{\tor}{{\rm Tor}\,}
\newcommand{\bff}{{\bf F}}
\newcommand{\bs}{\mathbb{S}}
\def\ns{\mathop{\lr}}
\def\nssup{\mathop{\lr\,sup}}
\def\nsinf{\mathop{\lr\,inf}}
\renewcommand{\phi}{\varphi}
\newcommand{\tT}{{\widetilde{T}}}
\newcommand{\tG}{{\widetilde{G}}}
\newcommand{\tB}{{\widetilde{B}}}
\newcommand{\tC}{{\widetilde{C}}}
\newcommand{\tW}{{\widetilde{W}}}
\newcommand{\tphi}{{\widetilde{\Phi}}}
\noindent

\section{Introduction}\label{Introduction}

We shall consider varieties over the field of complex numbers unless
otherwise specified.

Let $T\simeq (\mathbb{C}^*)^n$ denote the complex algebraic torus of
dimension $n$.  A {\em $T$-cellular variety} is a $T$-variety $X$
equipped with a $T$-stable algebraic cell decomposition. In other
words there is a filtration \be\label{filter} X=Z_1\supseteq
Z_2\supseteq\cdots\supseteq Z_m\supseteq Z_{m+1}=\emptyset\ee where
each $Z_i$ is a closed $T$-stable subvariety of $X$ and
$Z_i\setminus Z_{i+1}=Y_i$ is $T$-equivariantly isomorphic to the
affine space ${\mathbb{C}}^{n_i}$ equipped with a linear action of $T$
for $1\leq i\leq m$. Furthermore, $Y_i$ for $1\leq i\leq m$, are the
Bialynicki Birula cells associated to a {\em generic} one-parameter
subgroup of $T$ (see Definition \ref{cellular} and \cite{u4}).

We note here that any smooth projective complex variety $X$ with
$T$-action having only finitely many $T$-fixed points is $T$-cellular.
Indeed for any projective variety the Bialynicki-Birula cell
decomposition is filtrable (see \cite{BB1}) and since $X$ is smooth
the Bialynicki-Birula cells are smooth (see \cite[Section 3.1]{Br2}).

Let $T_{comp}\simeq (S^1)^n$ denote the maximal compact subgroup of
$T$.  Let $K^0$ denote the topological $K$-ring (see \cite{at}). For a
compact connected Lie group $G_{comp}$, let $K^0_{G_{comp}}$ denote
the $G_{comp}$-equivariant $K$-ring. In particular, $K^0_{T_{comp}}$
denotes the $T_{comp}$-equivariant topological $K$-ring (see
\cite{segal}).

We develop below some further notations to describe the results in
this paper.

Let $\mathtt{p}:\mathcal{E}\ra \mathcal{B}$ be a principal $T$-bundle,
where $T$ acts on $\mathcal{E}$ from the right. Let $G$ be a connected
complex reductive linear algebraic group which acts on $\mathcal{E}$ and
$\mathcal{B}$ from the left such that the bundle projection
$\mathtt{p}$ is $G$-equivariant. We note here that since the
$G$-action on $\mathcal{E}$ is from the left and $T$-action is from
the right, the actions commute.

For a $T$-cellular variety $X$ we consider the associated bundle
\[\pi:\mathcal{E}(X):=\mathcal{E}\times_{T} X\ra \mathcal{B},\] which is an
$G$-equivariant bundle with fiber $X$. We call $\mathcal{E}(X)$ a {\it
  cellular bundle}. Here $\pi([e,x])=\mathtt{p}(e)$ for
$[e,x]\in \mathcal{E}(X)$ (also see \cite{su}). We remark here that
a priori we do not assume the existance of $G$-action on the variety
$X$ and note that the action of $G$ on $\mathcal{E}(X)$ comes from the
left $G$-action on $\mathcal{E}$.

In \cite{u4}, we study the topological equivariant $K$-ring for
certain classes of singular toric varieties which are $T$-cellular
(see Definition \ref{cellular}). Let $X=X(\Delta)$ be a complete
$T$-cellular toric variety associated to the fan $\Delta$. There is a
combinatorial characterization on $\Delta$ so that $X(\Delta)$ is
$T$-cellular (see \cite[Theorem 3.1]{u4}). We consider the natural
restricted action of $T_{comp}\subset T$ on $X(\Delta)$.

In \cite[Theorem 4.2]{u4} (see Theorem \ref{main1} below), we
first give a GKM type description of the $T_{comp}$-equivariant
topological $K$-ring of $X(\Delta)$.  In \cite[Theorem 4.6]{u4}
(see Theorem \ref{main2} below), we further show that
$K^0_{T_{comp}}(X)$ is isomorphic as an $R(T_{comp})=R(T)$-algebra to
the ring of piecewise Laurent polynomial functions on $\Delta$ which
we denote by $PLP(\Delta)$.

\subsection{Our main results}

Let $X=X(\Delta)$ be a complete $T$-cellular toric variety. We shall
consider the $G$-equivariant associated toric bundle
$\mathcal{E}(X):=\mathcal{E}\times_{T} X$ over $\mathcal{B}$ with
bundle projection $\pi:\mathcal{E}(X)\ra \mathcal{B}$. 

Let $G_{comp}$ be a maximal compact subgroup of $G$.  

Our first aim in this paper is to extend the results in \cite{u4} to a
relative setting. More precisely, we describe the
$G_{comp}$-equivariant topological $K$-ring of
$\mathcal{E}(X(\Delta))$ as an algebra over the $G_{comp}$-equivariant
topological $K$-ring of $\mathcal{B}$.



Recall that $K^0_{G_{comp}}(\mathcal{B})$ gets a canonical
$R(T_{comp})$-algebra structure via the map
$R(T_{comp})\ra K^0_{G_{comp}}(\mathcal{B})$ which takes $e^{\chi}$ to
$[\mathcal{L}_{\chi}]_{G_{comp}}$, where $\mathcal{L}_{\chi}$ denotes
the $G_{comp}$-equivariant line bundle on $\mathcal{B}$
$\mathcal{E}\times_{T} \mathbb{C}_{\chi}$ for $\chi\in X^*(T_{comp})$.
Here $\mathbb{C}_{\chi}$ denotes the $1$-dimensional complex
representation of $T_{comp}$ associated to the character $\chi$ and
$T$ acts on $\mathbb{C}_{\chi}$ via its canonical projection to
$T_{comp}$.

Let $\mathcal{A}_{K^0_{G_{comp}}(\mathcal{B})}$ be the set of all
$(a_i)_{1\leq i\leq m}\in \Big(K^0_{G_{comp}}(\mathcal{B})\Big)^m$
such that $a_i\equiv a_j \pmod{(1-[\mathcal{L}_{\chi}]_{G_{comp}})}$
whenever the maximal dimensional cones $\sigma_i$ and $\sigma_j$ share
a wall $\sigma_i\cap \sigma_j$ in $\Delta$ and
$\chi\in (\sigma_i\cap \sigma_j)^{\perp}\cap M$. Moreover,
$\mathcal{A}_{K^0_{G_{comp}}(\mathcal{B})}$ is a
$K^0_{G_{comp}}(\mathcal{B})$-algebra where
$K^0_{G_{comp}}(\mathcal{B})$ is identified with the subalgebra of
$\mathcal{A}_{K^0_{G_{comp}}(\mathcal{B})}$ consisting of the diagonal
elements $(a,a,\ldots, a)$.

The following theorem gives a GKM type description of
  ${K}_{G_{comp}}^0(\mathcal{E}(X))$ as an
  ${K}_{G_{comp}}^0(\mathcal{B})$-algebra. This is a relative version
  of \cite[Theorem 4.2]{u4} (see Theorem \ref{main1} below).

  \bth (see Theorem \ref{main3}) Let $X=X(\Delta)$ be a complete
  $T$-cellular toric variety.  The ring
  $K^0_{G_{comp}}(\mathcal{E}(X))$ is isomorphic to
  $\mathcal{A}_{K^0_{G_{comp}}(\mathcal{B})}$ as an
  ${K}_{G_{comp}}^0(\mathcal{B})$-subalgebra of
  $\Big(K_{G_{comp}}^0(\mathcal{B})\Big)^m$. This isomorphism is
  obtained via the restriction \eqref{buninclusion} to the components
  $\mathcal{E}(x_i)$ for $1\leq i\leq m$ corresponding to the
  $T$-fixed points of $X$ (see Prop.
  \ref{relativelocalization}). \eeth


  Let $PLP(\Delta)$ be the ring of piecewise Laurent polynomial
  functions on $\Delta$ (see Section \ref{piecewiseLaurent} and
  \eqref{plp}).


  We have the following theorem which is a relative version of
  \cite[Theorem 4.6]{u4} (see Theorem \ref{main2} below).

  \bth (see Theorem \ref{cellulartorickunnethtop}) The ring
  ${K}_{G_{comp}}^0(\mathcal{E}(X))$ is isomorphic to
  ${K}_{G_{comp}}^0(\mathcal{B})\otimes_{R(T_{comp})}PLP(\Delta)$ as a
  ${K}_{G_{comp}}^0(\mathcal{B})$-algebra. This isomorphism is
  obtained via the restriction \eqref{buninclusion} to the components
  $\mathcal{E}(x_i)$ for $1\leq i\leq m$ corresponding to the
  $T$-fixed points of $X$ (see Prop.
  \ref{relativelocalization}). \eeth

  Taking $G=\{1\}$, we derive the ordinary topological $K$-ring of a
  cellular toric bundle in Corollary
  \ref{cellulartorickunnethtoptrivialH}. The results in Theorem
  \ref{cellulartorickunnethtop} and Corollary
  \ref{cellulartorickunnethtoptrivialH} generalize the results in
  \cite[Theorem 1.2 (ii), (iv)]{su} and \cite[Proposition 4.1]{u3} on
  the topological $K$-ring and Grothendieck ring of toric bundles with
  fibre a smooth projective $T$-toric variety to the setting of
  equivariant and ordinary topological $K$-ring of toric bundles with
  fibre any complete $T$-cellular toric variety. Furthermore, in
  Corollary \ref{eqkthtbsmooth}, taking $X$ to be smooth, we give a
  presentation for $K^0_{G_{comp}}(\mathcal{E}(X))$ as an algebra over
  $K^0_{G_{comp}}(\mathcal{B})$. This is the equivariant analogue of
  the results of Sankaran and the author \cite[Theorem 1.2(ii)]{su}.
  In the setting of non-smooth fibres the results on the description
  of the cohomology ring of weighted projective bundles are due to Al
  Amrani \cite{A3}. Also see the results by Bahri, Franz and Ray on
  the equivariant cohomology ring of weighted projective spaces and
  the cohomology ring of weighted projective bundles \cite{bfr}.

  For proving our first main results Theorem \ref{main3} and Theorem
  \ref{cellulartorickunnethtop}, our primary tool is a general result
  (see Theorem \ref{cellularkunnethtop} below), which is a Kunneth
  type formula which describes the $G_{comp}$-equivariant $K$-ring of
  a cellular bundle $\mathcal{E}(X)$ with fiber any $T$-cellular
  variety, as a tensor product of the $G_{comp}$-equivariant $K$-ring
  of the base $\mathcal{B}$ and the $T_{comp}$-equivariant $K$-ring of
  $X$, both of which have a canonical $R(T_{comp})$-algebra structure.
  We wish to mention that Theorem \ref{cellularkunnethtop} was proved
  in the non-equivariant setting i.e. when $G=1$ in \cite[Theorem
  2.4]{pu}. Also, the statement of the corresponding result for the
  equivariant setting was given in \cite[Theorem 3.2]{pu}, without
  proof. The above results in the non-equivariant and equivariant
  settings were respectively applied in \cite{pu}, by Paul and the
  author for the description of the ordinary and equivariant
  topological $K$-ring of Flag Bott manifolds of general Lie type.

Our second main aim in this paper is an application of the results on
the equivariant topological $K$-ring of cellular toric bundles proved
in Section \ref{mainresults}, to obtain a description of the
equivariant topological $K$-ring of toroidal horospherical
embedding. We outline these results below.

In Section \ref{horospherical}, we let $G$ to be a connected complex
reductive algebraic group and $H$ a closed subgroup of $G$. We say
that the closed subgroup $H$ is {\em horospherical} if it contains a
maximal unipotent subgroup $U$ of $G$. In this case we say that the
homogeneous space $G/H$ is {\em horospherical} (see
\cite{pasquier}). The normalizer $N_{G}(H)$ of a horospherical
subgroup $H$ in $G$ is a parabolic subgroup and $N_{G}(H)/H$ is a
torus (see \cite[Prop. 3.3]{monahan}). We call $P:=N_{G}(H)$ the {\em
  associated parabolic subgroup}, the flag variety $G/P$ the {\em
  associated flag variety} and $P/H$ the {\em associated torus} of the
horospherical homogeneous space $G/H$. Furthermore, the natural
projection map $G/H\ra G/P$ is a principal $P/H$-bundle.

Let $X$ be an equivariant embedding of the horospherical homogeneous
space $G/H$. Then $X$ is called a {\em horospherical variety} (see
\cite{knop}, \cite{t}, \cite{pasquier}, \cite{perrin} and
\cite{monahan}).

Recall that any toroidal horospherical embedding $X$ is of the form
$G\times_{P} Y$, which is a toric bundle with fibre a $P/H$-toric
variety $Y$, over the associated flag variety $G/P$ (see \cite[Chapter
13]{perrin}). Let $\mathbb{F}(X)$ denote the fan of $Y$ in the lattice
$N=X^*(P/H)$.

We shall further assume that $Y$ is $T$-cellular for the torus
$T=P/H$.

Using our description of the topological equivariant $K$-ring of a
cellular toric bundle in Theorem \ref{cellulartorickunnethtop}, in our
next main result Theorem \ref{equivhorospherical}, we describe the
$G_{comp}$-equivariant topological $K$-ring of the toroidal
horospherical variety $X$.

\bth(see Theorem \ref{equivhorospherical}) The $G_{comp}$-equivariant
topological $K$-ring of a complete toroidal horospherical embedding
$X$ of $G/H$ is isomorphic as an $K^0_{G_{comp}}(G/P)$-algebra to the
ring \[K^0_{G_{comp}}(G/P)\otimes_{R((P/H)_{comp})} PLP(\mathbb{F}(X)).\]\eeth

Furthermore, in the case when the toroidal horospherical variety $X$
is smooth which is equivalent to the condition that the $P/H$-toric
variety $Y$ is smooth, we get a presentation of the
$G_{comp}$-equivariant topological $K$-ring as well as the
$G$-equivariant Grothendieck ring of $X$, in Corollary
\ref{smoothcase} and Corollary \ref{smoothcaseGroth}
respectively. These are the $K$-theoretic analogues of the results of
Hofscheier, Khovanskii and Monin in \cite[Section 3.2]{kh} on the
cohomology ring of horospherical varieties.

\section{Preliminaries on $K$-theory}

Let $X$ be a compact $G_{comp}$-space for a compact Lie group
$G_{comp}$. By $K^0_{G_{comp}}(X)$ we mean the Grothendieck ring of
$G_{comp}$-equivariant topological vector bundles on $X$ with the
abelian group structure given by the direct sum and the multiplication
given by the tensor product of equivariant vector bundles. In
particular, $K^0_{G_{comp}}(pt)$, where $G_{comp}$ acts trivially on
$pt$, is the Grothendieck ring $R(G_{comp})$ of complex
representations of $G_{comp}$. The ring $K^0_{G_{comp}}(X)$ has the
structure of $R(G_{comp})$-algebra via the map
$R(G_{comp})\ra K_{G_{comp}}^0(X)$ which takes
$[V]\mapsto \mathbf{V}$, where $\mathbf{V}=X\times V$ is the trivial
$G_{comp}$-equivariant vector bundle on $X$ corresponding to the
$G_{comp}$-representation $V$. Let $pt$ be a $G_{comp}$-fixed point of
$X$ then the reduced equivariant $K$-ring
${\widetilde K}_{G_{comp}}^0(X)$ is the kernel of the map
$K_{G_{comp}}^0(X)\ra K_{G_{comp}}^0(pt)$, induced by the restriction
of $G_{comp}$-equivariant vector bundles. For $n\in \mathbb{N}$, we
define
$\widetilde{K}^{-n}_{G_{comp}}(X):=\widetilde{K}^0_{G_{comp}}(S^{n}X)$
where $S^{n}X$ is the $n$-fold reduced suspension of $X$. 

If $X$ is locally compact space but not compact we shall consider
$G_{comp}$-equivariant $K$-theory with compact support denoted
$K_{G_{comp},c}^0(X)$. This can be identified with
${\widetilde K}_{G_{comp}}^0(X^+)$ where $X^+$ denotes the one point
compactification of $X$, which is a compact $G_{comp}$-space, where
the point at infinity, which is the base point of $X^+$ is
$G_{comp}$-fixed. We define
${K}_{G_{comp},c}^{-n}(X):=\widetilde{K}^{-n}_{G_{comp}}(X^+)$.


When $X$ is already compact, we define $X^+=X\sqcup {pt}$ as the disjoint
union of $X$ and a base point. In this case we see that
$K^{0}_{G_{comp},c}(X)=\widetilde{K}^0_{G_{comp}}(X^+)=K^{0}_{G_{comp}}(X)$
since
$K^0_{G_{comp}}(X^+)=K_{G_{comp}}^0(X)\oplus K_{G_{comp}}^0(pt)$. In
\cite[Definition 2.8]{segal}, $K^{0}_{G_{comp},c}$ is denoted by
$K^{0}_{G_{comp}}$ without the subscript $c$. We remark here that
$K_{G_{comp},c}^0$ is not a homotopy invariant unlike $K_{G_{comp}}^0$
(see \cite[Proposition 2.3]{segal}). For example
$K^{0}_{G_{comp}, c}(\mathbb{R}^1)={\widetilde
  K}_{G_{comp}}^0(S^1)=K_{G_{comp}}^{-1}(pt)=0$ whereas
$K^0_{G_{comp},c}(pt)=K^0_{G_{comp}}(pt)=R(G_{comp})$.

We have the equivariant Bott periodicity
$K_{G_{comp}}^{-n}(X)\simeq K_{G_{comp}}^{-n-2}(X)$ given via
multiplication by the Bott element in $K_{G_{comp}}^{-2}(pt)$. (See
\cite{segal}, \cite{at} and \cite{fo}.) This enables us to define
$K^{n}_{G_{comp}}(X)$ for a positive $n\in \mathbb{Z}$ as
$K^{n-2q}(X)$ for $q\geq n/2$.

For $X$ a compact $G_{comp}$-space and $A$ a closed
$G_{comp}$-subspace for $n\in \mathbb{N}$ we define
$\widetilde{K}_{G}^{-n}(X, A):=\widetilde{K}^0(S^n(X/A))$ (see
\cite[Definition 2.7]{segal}, \cite[Definition 2.2.2]{at}) where $X/A$
is the space obtained by collapsing $A$ to a point. There is a long
exact sequence of $G_{comp}$-equivariant $K$-groups infinite in both
directions given as follows:

{\small \begin{align}\label{l.e.s} \cdots\ra
  {K}_{G_{comp}}^{-n}(X, A)\ra
  {K}_{G_{comp}}^{-n}(X) &\ra
  K_{G_{comp}}^{-n}(A)\ra \nonumber\\ & \ra
  K_{G_{comp}}^{-n+1}(X, A)\ra
  \cdots \end{align}}

Let $X$ be an algebraic variety with the action of an algebraic group
$G$. Then $\mathcal{K}^0_{G}(X)$ denotes the Grothendieck ring of
equivariant algebraic vector bundles on $X$ and $\mathcal{K}_0^{G}(X)$
the Grothendieck group of equivariant coherent sheaves on $X$ (see
\cite{Th}, \cite{mer}). The natural map
$\mathcal{K}^0_{G}(X)\ra \mathcal{K}_0^G(X)$ obtained by sending a
class of a $G$-equivariant vector bundle $\mathcal{V}$ on $X$ to the
dual of its sheaf of local sections is an isomorphism when $X$ is
smooth, but not in general. Moreover, when $G$ is a complex reductive
algebraic group and $G_{comp}$ is a maximal compact subgroup of $G$,
then any complex algebraic $G$-variety $X$ is a $G_{comp}$-space. When
$X$ is a smooth compact complex algebraic $G$-variety, we have a natural map
$\mathcal{K}^G_{0}(X)\ra K^0_{G_{comp}}(X)$, obtained by first
identifying $\mathcal{K}^{G}_0(X)$ with $\mathcal{K}_{G}^0(X)$ and
then viewing an algebraic $G$-vector bundle as a topological
$G_{comp}$-vector bundle on $X$ (see \cite[Section 5.5.5]{cg}).

Recall that in \cite[Section 5.5]{cg} we have an alternate notion of
$G$-cellular variety $X$, where $X$ admits a decreasing filtration
\eqref{filter} by $G$-stable closed subvarieties $Z_i$ such that
$Y_i=Z_i\setminus Z_{i+1}$ are complex affine spaces equipped with a
linear $G$-action. However, $Y_i$ are not assumed to be the cells of a
Bialynicki-Birula decomposition of $X$. Thus for $T$-varieties, this
notion of cellular is weaker than our definition (see Definition
\ref{cellular}). It follows from \cite[p. 272]{cg}, that the map
$\mathcal{K}^G_{0}(X)\ra K^0_{G_{comp}}(X)$ is an isomorphism when $X$
is smooth and $G$-cellular.

For singular $G$-varieties there are no natural isomorphisms
$\mathcal{K}^0_{G}(X)\cong \mathcal{K}_0^G(X)$ and
$\mathcal{K}^G_{0}(X)\cong K^0_{G_{comp}}(X)$.

\section{Bialynicki-Birula decomposition and cellular 
  varieties}\label{cellular varieties}

In this section we recall the notions of a filtrable Bialynicki-Birula
decomposition and a $T$-cellular variety (see \cite[Section 2]{u4} and
\cite[Section 3]{Br2}).

Let $X$ be a normal complex algebraic variety with the action of the
complex algebraic torus $T$. We assume that the set of $T$-fixed
points $X^{T}$ is finite. Let $\lambda$ be a generic one-parameter
subgroup of $T$ i.e.  $X^{\lambda}=X^{T}=\{ x_1,\ldots, x_m\}$. Let
\be\label{ps}\displaystyle X_{+}(x_i, \lambda)=\{ x\in X\mid
\lim_{t\ra 0}\lambda(t)x ~ \mbox{exists and is equal to}~x_i\}.\ee
Then $X_{+}(x_i,\lambda)$ is a locally closed $T$-invariant subvariety
of $X$ and is called the Bialynicki-Birula {\it plus} cell.


\bdefe\label{filtrable} The $T$-variety $X$ is called {\em
  filtrable} if it satisfies the following conditions:

(i) $X$ is the union of its plus cells $X_{+}(x_i,\lambda)$ for
$1\leq i\leq m$ with respect to a fixed generic one parameter subgroup
$\lambda$.

(ii) There exists a finite decreasing sequence of $T$-invariant closed
subvarieties of $X$
$X=Z_1\supset Z_2 \cdots\supset  Z_{m} \supset Z_{m+1}= \emptyset $
such that $Z_i\setminus Z_{i+1}=Y_i:=X_{+}(x_i,\lambda)$ for
$1\leq i\leq m$. In particular,
$\displaystyle \overline{Y_i}\subseteq Z_{i}=\bigcup_{j\geq i} Y_j$.
\edefe

\bdefe\label{cellular} Let $X$ be a normal complex algebraic variety
with an action of a complex algebraic torus $T$ such that
$X^T:=\{x_1,\ldots, x_m\}$ is finite. We say that $X$ is {\em $T$-cellular}
if $X$ is filtrable with respect to some fixed generic one-parameter
subgroup $\lambda$ and in addition each plus cell
$Y_i:=X_+(x_i,\lambda)$ is $T$-equivariantly isomorphic to a complex
affine space $\mathbb{C}^{n_i}$ on which $T$-acts linearly for
$1\leq i\leq m$.\edefe

\section{Equivariant $K$-theory of cellular
  bundles}\label{prelimkthcellularbundles}

We shall follow the notations in Section \ref{Introduction}.

Let $T\simeq (\mathbb{C}^*)^n$ and
$T_{comp}\simeq (S^1)^n\subseteq T$.  Let $G$ be a connected complex
reductive linear algebraic group. Let $\mathcal{B}$ be a nonsingular
complex algebraic $G$-variety, where $G$ acts on $\mathcal{B}$ from
the left. Let $\mathcal{E}$ be a nonsingular complex algebraic variety
with a right action of $T$ as well as a left action of $G$ such that
\[\mathtt{p}:{\mathcal E}\lra {\mathcal{B}=\mathcal{E}/T}\] is a principal $T$-bundle and
the projection $\mathtt{p}$ is $G$-equivariant.

Let $X$ be a $T$-cellular variety (see Definition \ref{filtrable}). We
shall consider the associated bundle
$\mathcal{E}(X):=\mathcal{E}\times_{T} X$. Then $G$ acts on
$\mathcal{E}(X)$ from the left and
$\pi:\mathcal{E}(X)\lra \mathcal{B}$ is $G$-equivariant. Let
$G_{comp}$ denote a maximal compact subgroup of $G$. Moreover, $\mathtt{p}$
and $\pi$ are $G_{comp}$-equivariant for the restricted
$G_{comp}$-action on $\mathcal{E}$ and $\mathcal{B}$ from the left.

We have the following decreasing sequence of
$G$-equivariant bundles over $\mathcal{B}$: \be\label{relfilt}
\mathcal{E}(X)\supseteq \mathcal{E}(Z_1)\supseteq \cdots \supseteq
\mathcal{E}(Z_m) \supseteq \mathcal{E}(Z_{m+1})=
\emptyset\ee where
$\mathcal{E}(Z_i)\setminus \mathcal{E}(Z_{i+1})=\mathcal{E}(Y_i)$ is
an $G$-equivariant affine bundle over $\mathcal{B}$, associated to the
linear $T$-representation $\mathbb{C}^{k_i}$.

The $G$-equivariant principal $T$-bundle
$\mathtt{p}:\mathcal{E}\ra \mathcal{B}$ is Zariski locally trivial and
the transition functions of the associated affine bundle are also
affine.

\begin{assum}\label{base}
  We shall assume throughout this section that
  \[K_{G_{comp}}^{-1}(\mathcal{B})=0\] and $K^{-1}(\mathcal{B})=0$. We
  shall also assume that $X$ and $\mathcal{B}$ are compact topological
  spaces. In particular, this implies that for each $1\leq i\leq m$,
  $Z_i$ and $\mathcal{E}(Z_i)$ are compact. Moreover,
  $Z_m=Y_m=\{x_{m}\}$ and
  $\mathcal{E}(Z_m)=\mathcal{E}({x_m})=\mathcal{B}$. \end{assum}

The following result is a relative version of \cite[Theorem 3.1]{u4}
in the setting of $K$-theory. In the non-equivariant setting i.e. when
$G=1$, this result was proved in \cite[Theorem 2.4]{pu}, and for the
equivariant setting it was stated without proof (see \cite[Theorem
3.2]{pu}). However in \cite{pu}, the Assumption \ref{base} was
implicitly assumed and not explicilty stated before the statement of
\cite[Theorem 3.2]{pu}. We state and give a detailed proof of this
result to make the exposition self contained. We wish to emphasise
that the arguments used in the proof are essentially an
$G_{comp}$-equivariant analogue of those used in the proof of
\cite[Theorem 2.4]{pu}. We have shortened the arguments in some places
by giving references to \cite[Theorem 3.1]{u4} and \cite[proof of
Theorem 2.4]{pu}.

Consider the map \be\label{topkunneth}\varphi:
{K}_{G_{comp}}^0(\mathcal{B})\otimes_{R(T_{comp})}{K}^0_{T_{comp}}(X)\ra
{K}_{G_{comp}}^0(\mathcal{E}(X))\ee which maps
$[\mathcal{V}]\otimes [\mathcal{W}]$ to
$[\pi^*(\mathcal{V})]\cdot [\mathcal{E}\times_{T} \mathcal{W}]$ where
$[\mathcal{V}]$ is the class in ${K}_{G_{comp}}^0(\mathcal{B})$ of an
$G_{comp}$-equivariant vector bundle $\mathcal{V}$ on $\mathcal{B}$
and $[\mathcal{W}]$ is the class in ${K}^0_{T_{comp}}(X)$ of a
$T_{comp}$-equivariant vector bundle $\mathcal{W}$ on $X$. Note that
$T$ acts on $\mathcal{W}$ via its canonical projection to $T_{comp}$.

Note that ${K}_{G_{comp}}^0(\mathcal{B})$ is a $R(T_{comp})$-module
via the map which takes any $[\mathbf{V}]\in R(T_{comp})$ for a
$T_{comp}$-representation $\mathbf{V}$ to
$[\mathcal{E}\times_{T} \mathbf{V}]\in
{K}_{G_{comp}}^0(\mathcal{B})$. Here $T$ acts on $\mathbf{V}$ via its
projection to $T_{comp}$. Furthermore, $G_{comp}$ acts on
$\mathcal{E}(\mathbf{V}):=\mathcal{E}\times_{T_{comp}} \mathbf{V}$
from the left.

\bth\label{cellularkunnethtop} (\cite[Theorem 3.2]{pu}) The map
$\varphi$ is an isomorphism of
${K}_{G_{comp}}^0(\mathcal{B})$-algebras.  Further, we have
$K_{G_{comp}}^{-1}(\mathcal{E}(X))=0$.\eeth

\begin{proof}

  Following the arguments in the proof of \cite[Theorem 3.1]{u4}, when
  $q$ is even \[K^{-q}_{T_{comp}}(Z_i,Z_{i+1})=R(T_{comp})\] and when
  $q$ is odd \[K^{-q}_{T_{comp}}(Z_i,Z_{i+1})=0.\] In particular, we
  get the following split short exact sequence of
  $R(T_{comp})$-modules {\small \be \label{eq1} 0\ra
    K_{T_{comp}}^0(Z_i, Z_{i+1})\ra K_{T_{comp}}^0(Z_i)\stackrel
    {\alpha^*} {\ra} K_{T_{comp}}^0( Z_{i+1})\ra 0 \ee} for
  $1\leq i\leq m$.  Furthermore, it also follows from \eqref{eq1} and
  descending induction on $i$ that $K_{T_{comp}}^0(Z_i)$ is a free
  $R(T_{comp})$-module of rank $m-i+1$ and $K^{-1}_{T_{comp}}(Z_i)=0$
  for $1\leq i\leq m$. We can start the descending induction since
  $Z_m=Y_m=\{x_m\}$, so that $K_{T_{comp}}^0(Z_m)=R(T_{comp})$ and
  $K_{T_{comp}}^{-1}(Z_m)=0$.

Now, by tensoring \eqref{eq1} with $K_{G_{comp}}^0(\mathcal{B})$ over
$R(T_{comp})$ we get the exact sequence
{\tiny \be\label{eq2} \displaystyle 0\ra K_{G_{comp}}^0(\mathcal{B})\ra
K_{G_{comp}}^0(\mathcal{B}) \bigotimes_{R(T_{comp})}
K^0_{T_{comp}}(Z_i)\ra K_{G_{comp}}^0(\mathcal{B})
\bigotimes_{R(T_{comp})} K^0_{T_{comp}}(Z_{i+1})\ra 0\ee} for
$1\leq i\leq m$.

Since $\mathcal{E}(Z_{i+1})\subseteq \mathcal{E}(X_{i})$ is closed
subspace for $1\leq i\leq m$, by \cite[Proposition 2.4.4, Theorem
2.4.9]{at} and \cite{segal} it follows that we have the following long
exact sequence of $G_{comp}$-equivariant
$K$-groups: {\small \begin{align}\label{les1} \cdots\ra
  {K}_{G_{comp}}^{-q}(\mathcal{E}(Z_i), \mathcal{E}(Z_{i+1}))\ra
  {K}_{G_{comp}}^{-q}(\mathcal{E}(Z_i)) &\ra
  K_{G_{comp}}^{-q}(\mathcal{E}(Z_{i+1}))\ra \nonumber\\ & \ra
  K_{G_{comp}}^{-q+1}(\mathcal{E}(Z_i), \mathcal{E}(Z_{i+1}))\ra
  \cdots \end{align}} for $1\leq i\leq m$.

Furthermore, we have
\begin{align*} K_{G_{comp}}^{-q}(\mathcal{E}(Z_i),
  \mathcal{E}(Z_{i+1}))&=\widetilde{K}_{G_{comp}}^{-q}(\mathcal{E}(Z_i)/
  \mathcal{E}(Z_{i+1}))\\
  &=\widetilde{K}_{G_{comp}}^{-q}(\mathcal{E}(Y_i)^+)\\ &\simeq
  K_{G_{comp}}^{-q}(\mathcal{B})\end{align*} for $1\leq i\leq m$ (see
\cite[Definition 2.4.1, Proposition 2.7.2]{at}).

By descending induction on $i$ we shall assume that
$K_{G_{comp}}^{-1}(\mathcal{E}(Z_{i+1}))=0$ and
$K_{G_{comp}}^0(\mathcal{E}(Z_{i+1}))$ is a free
$K_{G_{comp}}^0(\mathcal{B})$-module of rank $m-i$. We can start the
induction since
\[K_{G_{comp}}^{-1}(\mathcal{E}(Z_m))=K_{G_{comp}}^{-1}(\mathcal{E}(Y_m))=K_{G_{comp}}^{-1}(\mathcal{B})=0\]
and
\[K_{G_{comp}}^0(\mathcal{E}(Z_m))=K_{G_{comp}}^0(\mathcal{E}(Y_m))=K_{G_{comp}}^0(\mathcal{B})\] (by Assumption \ref{base}). Then
from \eqref{les1} we get the following short exact sequence

{\small \be\label{eq3} 0\lra K_{G_{comp}}^0(\mathcal{B})\lra
K_{G_{comp}}^0(\mathcal{E}(Z_i))\lra
K_{G_{comp}}^0(\mathcal{E}(Z_{i+1}))\lra 0.\ee} By induction we also
get that $K_{G_{comp}}^{-1}(\mathcal{E}(Z_i))=0$, since
\[K_{G_{comp}}^{-1}(\mathcal{E}(Z_i),
  \mathcal{E}(Z_{i+1}))=\widetilde{K}_{G_{comp}}^{-1}(\mathcal{E}(Y_i)^+)=
  K_{G_{comp}}^{-1}(\mathcal{B})=0\]
(by Assumption \ref{base}).

Now, \eqref{eq3} splits as a
$K_{G_{comp}}^0(\mathcal{B})$-module. Hence by induction it follows
that $K_{G_{comp}}^0(\mathcal{E}(Z_i))$ is a free
$K_{G_{comp}}^0(\mathcal{B})$ module of rank $m-i+1$ and
$K_{G_{comp}}^{-1}(\mathcal{E}(Z_i))=0$.  In particular, when $i=1$,
it follows that $K_{G_{comp}}^0(\mathcal{E}(X))$ is a free
$K_{G_{comp}}^0(\mathcal{B})$-module of rank $m$ and
$K_{G_{comp}}^{-1}(\mathcal{E}(X))=0$.

We shall now show by descending induction on $i$ that the map
\[\phi_i: K_{G_{comp}}^0(\mathcal{B}) \otimes_{R(T_{comp})} K^0_{T_{comp}}(Z_i)\lra
  K_{G_{comp}}^0(\mathcal{E}(Z_i))\] defined by sending
$[\mathcal{V}]\otimes [\mathcal{W}]$ to
$[\pi^*(\mathcal{V})]\cdot [\mathcal{E}\times_{T} \mathcal{W}]$
is an isomorphism for $1\leq i\leq m$. Here $[\mathcal{W}]$ is the
class in $K^0_{T_{comp}}(Z_i)$ of a $T_{comp}$-equivariant vector
bundle on $Z_i$ and $[\mathcal{V}]$ is the class in
$K_{G_{comp}}^0(\mathcal{B})$ of an $G_{comp}$-equivariant vector
bundle on $\mathcal{B}$. This will prove the theorem since
$\phi_1=\phi$. We can start the induction since
\[\phi_m:K_{G_{comp}}^0(\mathcal{B})\bigotimes_{R(T_{comp})} K^0_{T_{comp}}(Z_m)\lra
  K_{G_{comp}}^0(\mathcal{E}(Y_m))=K_{G_{comp}}^0(\mathcal{B})\] is
an isomorphism. This follows because $Z_m=Y_m=\{x_m\}$ and
$K^0_{T_{comp}}(Y_m)=R(T_{comp})$. Consider the following commuting
diagram where the first row is \eqref{eq2} and the second row is
\eqref{eq3}.  {\tiny \begin{center}
	\begin{tikzcd}[cramped,column sep=tiny]
          0\ar[r] & K_{G_{comp}}^0(\mathcal{B}) \arrow[r]
          \arrow[d,-,double equal sign distance,double] &
          K_{G_{comp}}^0(\mathcal{B})\bigotimes_{R(T_{comp})}K^0_{T_{comp}}(Z_i)
          \arrow[r] \arrow[d, "\phi_i"]
          &K_{G_{comp}}^0(\mathcal{B})\bigotimes_{R(T_{comp})}K^0_{T_{comp}}(Z_{i+1})
          \arrow[d, "\phi_{i+1}"]\ar[r]
          &0\\
          0\ar[r] & K_{G_{comp}}^0(\mathcal{B})\arrow[r]
          &K_{G_{comp}}^0(\mathcal{E}(Z_i))\arrow[r]
          &K_{G_{comp}}^0(\mathcal{E}(Z_{i+1}))\ar[r]
          &0\\
	\end{tikzcd}
\end{center}}

In the above diagram the first vertical
arrow is the identity map and $\phi_{i+1}$ is an isomorphism by
induction hypothesis. By diagram chasing it follows that $\phi_i$ is
an isomorphism. Hence the theorem.

\end{proof}

In particular, when  $G=\{1\}$, the above theorem reduces to the following corollary.

\bcor\label{cellularkunnethtoptrivialH}(\cite[Theorem 2.4]{pu}) We have the
following isomorphism: \be\label{topkunnethtriv}{K}^0(\mathcal{E}(X))\simeq
{K}^0(\mathcal{B})\otimes_{R(T_{comp})}{K}^0_{T_{comp}}(X).\ee Here ${K}^0(\mathcal{B})$
is a $R(T_{comp})$-module via the map which takes any
$[\mathbf{V}]\in R(T_{comp})$ for a $T$-representation $\mathbf{V}$ to
$[\mathcal{E}\times_{T} \mathbf{V}]\in {K}^0(\mathcal{B})$.  Further, we have
$K^{-1}(\mathcal{E}(X))=0$. \ecor

Furthermore, when $\mathcal{B}=pt$ we have
$K^0(\mathcal{B})=\mathbb{Z}$, $K^{-1}(\mathcal{B})=0$ and
$\mathcal{E}(X)\cong X$. Thus the above theorem reduces to the
following corollary. In other words $T$-cellular varieties are {\em
  weakly equivariantly formal in $K$-theory} with respect to
$T_{comp}$-action, in the sense of \cite{hl}.

\bcor\label{tequivformal} We have the isomorphism \be\label{teqf}
\mathbb{Z}\otimes_{R(T_{comp})}{K}^0_{T_{comp}}(X) \cong {K}^0(X).\ee
Here the map from $K_{T_{comp}}(X)\lra K(X)$ is the canonical
forgetful map, $\mathbb{Z}\lra K^0(X)$ is the map which takes an
integer $n$ to the class of the trivial vector bundle of rank $n$ and
$\mathbb{Z}$ is $R(T_{comp})$-module under the augmentation map. \ecor

In the following proposition we isolate the important result for a
$T$-cellular variety $X$ (see Definition \ref{cellular}) which was
also proved during the course of the proof of Theorem
\ref{cellularkunnethtop}.

\bpropo\label{celldec}The ring $K^0_{T_{comp}}(X)$ is a free
$R(T_{comp})$-module of rank $m$ which is the number of
cells. Furthermore, we have $K^1_{T_{comp}}(X)=0$.  \epropo

We have the following corollary for any $T$-cellular variety $X$.

Let $\iota:X^{T_{comp}}=X^{T}\hra X$ denote the inclusion of the set of
$T$-fixed points in $X$.

\bcor\label{localization} (see \cite[Corollary 3.2]{u4}) The canonical
restriction map
\[ K^0_{T_{comp}}(X)\stackrel{\iota^*}{\lra} K^0_{T_{comp}}(X^{T_{comp}})\cong R(T_{comp})^m\] is
injective where $m=|X^{T_{comp}}|$.  \ecor

Let $X^{T}=\{x_1,\ldots, x_m\}$ denote the the $T$-fixed points of
  $X$. For $1\leq i\leq m$, the map $s_i: \mathcal{B}\ra \mathcal{E}(X)$
  defined by $b\mapsto [e,x_i]$ for $e\in \mathtt{p}^{-1}(b)$ defines
  a $G$-equivariant section for the bundle
  $\pi:\mathcal{E}(X)\ra \mathcal{B}$. In particular,
  $\pi\circ s_i=Id_{\mathcal{B}}$ so that
  $\pi^*\circ s_i^{*}=Id^*_{K^0_{H_{comp}}(\mathcal{B})}$. In
  particular,
  $\pi^*: K^0_{G_{comp}}(\mathcal{B})\ra
  K^0_{G_{comp}}(\mathcal{E}(X))$ is a monomorphism of
  $R(G_{comp})$-algebras. Thus we have a canonical
  $K^0_{G_{comp}}(\mathcal{B})$-algebra structure on the ring
  $K^0_{G_{comp}}(\mathcal{E}(X))$.

  The inclusion $X^{T}\hra X$ induces an $G$-equivariant inclusion
  \be\label{buninclusion}\iota_{rel, X}: \mathcal{E}(X^{T})\hra
  \mathcal{E}(X).\ee Furthermore, since
  $\displaystyle \mathcal{E}(X^{T})=\bigsqcup_{i=1}^m
  \mathcal{E}\times_{T} x_i$ and
  $\mathcal{E}(x_i):=\mathcal{E}\times_{T} x_i\cong \mathcal{B}$ for
  $1\leq i\leq m$, it follows that
  \[K^0_{G_{comp}}(\mathcal{E}(X^{T}))\cong
   \Big(K_{G_{comp}}^0(\mathcal{B})\Big)^m.\]

The following is a relative version of Corollary \ref{localization}.
  
\bpropo\label{relativelocalization} The map
$\iota_{rel, X}^*: K^0_{G_{comp}}(\mathcal{E}(X))\ra K^0_{G_{comp}}(\mathcal{E}(X^T))$ is an
injective map of $K^0_{G_{comp}}(\mathcal{B})$-algebras.
\epropo
 \begin{proof} 
   Consider the following commuting diagram where the first row is
   \eqref{eq3}, and the vertical arrows are the maps
   $ \iota^*_{rel, Z_i}$ and $\iota^*_{rel, Z_{i+1}}$, induced by
   restriction to the $T$-fixed points of $Z_i$ and $Z_{i+1}$
   respectively:

{\tiny \begin{center}
	\begin{tikzcd}[cramped,column sep=tiny]
          0\ar[r] & K_{G_{comp}}^0(\mathcal{B})\arrow[r]
          \arrow[d,-,double equal sign distance,double]
          &K_{G_{comp}}^0(\mathcal{E}(Z_i))\arrow[r] \arrow[d ]
          &K_{G_{comp}}^0(\mathcal{E}(Z_{i+1})) \arrow[d]\ar[r] &0\\
          0\ar[r] & K_{G_{comp}}^0(\mathcal{B}) \arrow[r] &
          K^0_{G_{comp}}(\mathcal{E}(Z_i^{T})) \arrow[r] &
          K^0_{G_{comp}}(\mathcal{E}(Z_{i+1}^{T}))
          \ar[r] &0 \\
	\end{tikzcd}
\end{center}}

The proof of injectivity of
$\iota^*_{rel, Z_i}: K_{G_{comp}}^0(\mathcal{E}(Z_i))\ra
K^0_{G_{comp}}(\mathcal{E}(Z_i)^{T})=
\Big(K_{G_{comp}}^0(\mathcal{B})\Big)^{m-i+1}$ follows by descending
induction on $i$ and diagram chasing.  We can start the induction
since $Z_m=\{x_m\}$ and
$\iota^*_{rel, Z_m}=id_{K_{G_{comp}}^0(\mathcal{B})}$. Since $Z_1=X$
the proposition follows.
\end{proof}

\section{GKM theory of  Cellular Toric varieties}\label{GKMtoricsection}




We begin by fixing some notations and conventions.

Let $X=X(\Delta)$ be the toric variety associated to a fan $\Delta$ in
the lattice $N\simeq \bz^n$. Let $M:=\mbox{Hom}(N,\bz)$ be the dual
lattice of characters of $T$. Let $\{v_1,\ldots,v_d\}$ denote the set
of primitive vectors along the edges
$\Delta(1):=\{\rho_1,\ldots,\rho_d\}$. Let $V(\gamma)$ denote the
orbit closure in $X$ of the $T$-orbit $O_{\gamma}$ corresponding to
the cone $\gamma\in\Delta$. Let $S_{\sigma}=\sigma^{\vee}\cap M$ and
$U_{\sigma}:=\mbox{Hom}_{sg}(S_{\sigma},\mathbb{C})$ denote the
$T$-stable open affine subvariety corresponding to a cone
$\sigma\in \Delta$. Here $\mbox{Hom}_{sg}$ denotes semigroup
homomorphisms.

We further assume that all the maximal cones in $\Delta$ are
$n$-dimensional, in other words $\Delta$ is {\it pure}. The
$T$-fixed locus in $X$ consists of the set of $T$-fixed points
\be\label{fp}\{x_{1}, x_{2}\ldots,x_{m}\}\ee
corresponding to the set of maximal dimensional cones
\be\label{maxc}\Delta(n):=\{\sigma_1,\sigma_2,\ldots,\sigma_m\}.\ee

In \cite[Theorem 3.1]{u4} a combinatorial characterization is given on
the fan $\Delta$ for the toric variety $X(\Delta)$ to be $T$-cellular
(see Definition \ref{cellular}). We call a fan $\Delta$ satisfying the
necessary and sufficient combinatorial conditions a {\em cellular fan}
(see \cite[Definition 3.2]{u4}).

Henceforth we assume that $\Delta$ is a complete cellular fan or
equivalently that $X:=X(\Delta)$ is a complete $T$-cellular toric
variety.


Recall that (see \cite[Section 5]{u4} and \cite[272-273]{cg}) we have
the isomorphism \be\label{iso3'} R(T)\cong R(T_{comp})=\mathbb{Z}[M]
\ee where $\mathbb{Z}[M]=\mathbb{Z}[e^{u} : u\in M]$.

Let $\mathcal{A}$ be the set of all
$(a_i)_{1\leq i\leq m}\in R(T_{comp})^m$ such that
$a_i\equiv a_j \pmod{(1-e^{\chi})}$ whenever the maximal dimensional
cones $\sigma_i$ and $\sigma_j$ share a wall $\sigma_i\cap \sigma_j$
in $\Delta$ and $\chi\in (\sigma_i\cap \sigma_j)^{\perp}\cap M$. In
other words, the $T$-fixed points $x_i$ and $x_j$ lie in the closed
$T$-stable irreducible curve $C_{ij}:=V(\sigma_i\cap \sigma_j)$ in
$X(\Delta)$ and $T$ (and hence $T_{comp}$) acts on $C_{ij}$ through
the character $\chi$. Moreover, $\mathcal{A}$ is an
$R(T_{comp})$-algebra where $R(T_{comp})$ is identified with the
subalgebra of $\mathcal{A}$ consisting of the diagonal elements
$(a,a,\ldots, a)$.

The following theorem gives a GKM type description for
$K^0_{T_{comp}}(X)$ as a $K^0_{T_{comp}}(pt)=R(T_{comp})$-subalgebra
of $R(T_{comp})^m$ (see Cor. \ref{localization}).

\bth\label{main1}(\cite[Theorem 4.2]{u4}) Let $X=X(\Delta)$ be a
complete $T$-cellular toric variety.  The ring $K^0_{T_{comp}}(X)$ is
isomorphic to $\mathcal{A}$ as an $R(T_{comp})$-subalgebra of
$R(T_{comp})^m$.  \eeth


\subsection{Piecewise Laurent polynomial functions on
  $\Delta$}\label{piecewiseLaurent}

Let $X=X(\Delta)$ be a complete $T$-cellular toric variety.






Let $\sigma\in \Delta$.  Then $\sigma^{\perp}\cap M$ is a sublattice
of $M$. We have the identifications \be\label{ident}
K_{T_{comp}}^0(O_{\sigma})=K_{T_{comp}}^0(T/T_{\sigma})=R(T_{\sigma})\ee
for $\sigma\in \Delta$ where $O_{\sigma}$ is the $T$-orbit in $X$
corresponding to $\sigma$ and $T_{\sigma}\subseteq T$ is the
stabilizer of $O_{\sigma}$.

Further, since $X^*(T_{\sigma})=M/M\cap \sigma^{\perp}$ we can
identify $R(T_{\sigma})$ with $\mathbb{Z}[M/M\cap \sigma^{\perp}]$
which is the ring of Laurent polynomial functions on $\sigma$, for
$\sigma\in \Delta$.  Furthermore, whenever $\sigma$ is a face of
$\sigma'\in\Delta$ we have a natural homomorphism
\be\label{restLaurent} \psi_{\sigma,\sigma'}:\mathbb{Z}[M/\sigma'^{\perp}\cap
M]=R(T_{\sigma'})\ra \mathbb{Z}[M/\sigma^{\perp}\cap
M]=R(T_{\sigma})\ee given by the restriction of Laurent polynomial
functions on $\sigma'$ to $\sigma$.

Let \be\label{plp}PLP(\Delta):=\{(f_{\sigma}) \in \prod_{\sigma\in
  \Delta} {\mathbb{Z}[M/M\cap
  \sigma^{\perp}]}~\mid~\psi_{\sigma,\sigma'} (f_{\sigma'})=
f_{\sigma}~ \mbox{whenever}~ \sigma\preceq \sigma' \in \Delta\}.\ee
Then $PLP(\Delta)$ is a ring under pointwise addition and
multiplication and is called the ring of piecewise Laurent polynomial
functions on $\Delta$. Moreover, we have a canonical map
$R(T_{comp})=\mathbb{Z}[M]\ra PLP(\Delta)$ which sends $f$ to the
constant tuple $(f)_{\sigma\in \Delta}$. This gives $PLP(\Delta)$ the
structure of $R(T_{comp})$ algebra. We have the following theorem.

\bth\label{main2} (see \cite[Theorem 4.6]{u4})The ring
$K^0_{T_{comp}}(X(\Delta))$ is isomorphic to $PLP(\Delta)$ as an
$R(T_{comp})$-algebra.  \eeth
  
\brem The description of $K^0_{T_{comp}}(X)$ in Theorem \ref{main1}
and Theorem \ref{main2} agrees with the results of Anderson and Payne
(see \cite[Theorem 1.6]{ap}) on the algebraic equivariant operational
$K$-ring of $X(\Delta)$. Further, these extend the results for the
topological equivariant $K$-ring of a retractable toric orbifolds
(which include cellular projective simplicial toric varieties) due to
Sarkar and the author in \cite[Theorem 4.2]{saruma} and the
topological equivariant $K$-ring of a divisive weighted projective
space due to Harada, Holm, Ray and Williams in \cite[Theorem
5.5]{hhrw}. In the non-equivariant setting, the results on the
$K$-theory of weighted projective spaces are due to Al Amrani (see
\cite{A1, A2}). The results on the algebraic equivariant $K$-ring of
smooth toric varieties are due to Vezzosi and Vistoli \cite{vv}. The
corresponding results on the topological $T_{comp}$-equivariant
$K$-ring of a smooth projective toric variety were proved in
\cite[Theorem 7.1, Corollary 7.2]{hhrw}.  In the non-equivariant
setting Sankaran and the author in \cite{su} gave presentations for
the topological $K$-ring as well as the Grothendieck ring of smooth
projective toric varieties in terms of generators and relations (also
see \cite{s} for the corresponding results on smooth complete toric
varieties).  \erem

\section{Equivariant $K$-theory of cellular toric bundles}\label{mainresults}

Following the notations of Section \ref{prelimkthcellularbundles}, $H$
is a connected linear algebraic group, $\mathcal{E}$ is a nonsingular
complex algebraic variety with a left action of $H$ and a right action
of $T$ such that
  \[\mathtt{p}:{\mathcal E}\lra {\mathcal{B}=\mathcal{E}/T}\] be an
  $G$-equivariant principal $T$-bundle, where $\mathcal{B}$ is a
  nonsingular algebraic variety with a left $G$-action satisfying
  Assumption \ref{base}.

  Let $\mathcal{L}_{\chi}$ denote the $G_{comp}$-equivariant line
  bundle $\mathcal{E}\times_{T_{comp}} \mathbb{C}_{\chi}$ for
  $\chi\in X^*(T)$.

  Recall that $K^0_{G_{comp}}(\mathcal{B})$ gets a $R(T_{comp})$-algebra
  structure via the map $R(T_{comp})\ra K^0_{G_{comp}}(\mathcal{B})$
  which takes $e^{\chi}$ to $[\mathcal{L}_{\chi}]_{G_{comp}}$.
  
  Let $X=X(\Delta)$ be a $T$-cellular toric variety. We shall consider
  the $H$-equivariant associated toric bundle
  $\mathcal{E}(X):=\mathcal{E}\times_{T} X$ over $\mathcal{B}$ with
  bundle projection $\pi:\mathcal{E}(X)\ra \mathcal{B}$.

  Let $\mathcal{A}_{K^0_{G_{comp}}(\mathcal{B})}$ be the set of all
  $(a_i)_{1\leq i\leq m}\in K^0_{G_{comp}}(\mathcal{B})^m$ such that
  $a_i\equiv a_j \pmod{(1-[\mathcal{L}_{\chi}]_{H_{comp}})}$ whenever
  the maximal dimensional cones $\sigma_i$ and $\sigma_j$ share a wall
  $\sigma_i\cap \sigma_j$ in $\Delta$ and
  $\chi\in (\sigma_i\cap \sigma_j)^{\perp}\cap M$. In other words, the
  $T$-fixed points $x_i$ and $x_j$ lie in the closed $T$-stable
  irreducible curve $C_{ij}:=V(\sigma_i\cap \sigma_j)$ in $X(\Delta)$
  and $T$ (and hence $T_{comp}$) acts on $C_{ij}$ through the
  character $\chi$. Moreover,
  $\mathcal{A}_{K^0_{G_{comp}}(\mathcal{B})}$ is a
  $K^0_{G_{comp}}(\mathcal{B})$-algebra where
  $K^0_{G_{comp}}(\mathcal{B})$ is identified with the subalgebra of
  $\mathcal{A}_{K^0_{G_{comp}}(\mathcal{B})}$ consisting of the
  diagonal elements $(a,a,\ldots, a)$.

  The following theorem gives a GKM type description of
  ${K}_{G_{comp}}^0(\mathcal{E}(X))$ as an
  ${K}_{G_{comp}}^0(\mathcal{B})$-algebra.

  \bth\label{main3} Let $X=X(\Delta)$ be a complete $T$-cellular toric
  variety.  The ring $K^0_{G_{comp}}(\mathcal{E}(X))$ is isomorphic to
  $\mathcal{A}_{K^0_{G_{comp}}(\mathcal{B})}$ as an
  ${K}_{G_{comp}}^0(\mathcal{B})$-algebra.  The isomorphism is given
  via restriction \eqref{buninclusion} to the components
  $\mathcal{E}(x_i)$ for $1\leq i\leq m$ corresponding to the
  $T$-fixed points of $X$ (see Prop.
  \ref{relativelocalization}).\eeth
  \begin{proof} The map
    \[\iota^*_{res, X}: K^0_{G_{comp}}(\mathcal{E}(X))\ra
      K^0_{G_{comp}}(\mathcal{E}(X^{T}))\cong
      \Big(K^0_{G_{comp}}(\mathcal{B})\Big)^{m}\] is an inclusion by
    Prop. \ref{relativelocalization}. Moreover, by Theorem
    \ref{cellularkunnethtop}, it follows that $\iota^*_{res, X}$ is
    induced from
    $\iota^*: K^0_{T_{comp}}(X)\hra K^0_{T_{comp}}(X^{T})\cong
    \Big(R(T_{comp})\Big)^{m}$, by extending scalars to the
    $R(T_{comp})$-algebra $K^0_{G_{comp}}(\mathcal{B})$. The theorem
    now follows from Theorem \ref{main1} and the observation that
    \[\mathcal{A}_{K^0_{G_{comp}}(\mathcal{B})}\cong
      K^0_{G_{comp}}(\mathcal{B})\bigotimes_{R(T_{comp})}\mathcal{A}.\]
  \end{proof}


  \brem\label{comphaheho} The above result can alternately be proved
  by using the theorem \cite[Theorem 3.1]{haheho} of Harada Henriques
  and Holm. Our proof is based on the Kunneth type formula namely
  Theorem \ref{cellularkunnethtop}, and Theorem \ref{main1} which was
  proved in \cite[Theorem 4.2]{u4}. For an algebraic $K$-theoretic
  version of this result we refer to \cite[Theorem 1.11, Section
  1.2]{u3} \erem

The following theorem describes ${K}_{G_{comp}}^0(\mathcal{E}(X))$ as
a piecewise ${K}_{G_{comp}}^0(\mathcal{B})$-algebra over the fan
$\Delta$.

\bth \label{cellulartorickunnethtop} Let $X=X(\Delta)$ be a complete
$T$-cellular toric variety. We have the following isomorphism of
${K}_{G_{comp}}^0(\mathcal{B})$-algebras: \be
{K}_{G_{comp}}^0(\mathcal{E}(X))\simeq
{K}_{G_{comp}}^0(\mathcal{B})\otimes_{R(T_{comp})}PLP(\Delta).\ee This isomorphism is
  obtained via the restriction \eqref{buninclusion} to the components
  $\mathcal{E}(x_i)$ for $1\leq i\leq m$ corresponding to the
  $T$-fixed points of $X$ (see Prop.
  \ref{relativelocalization}). 
\eeth

  \begin{proof}
    This follows from Theorem \ref{cellularkunnethtop}, using the
    description of $K^0_{T_{comp}}(X(\Delta))$ given in Theorem
    \ref{main2}.

  \end{proof}

  Further, putting $G=\{1\}$, and using Corollary
  \ref{cellularkunnethtoptrivialH}, we get the following corollary
  which is a generalization of \cite[Theorem 1.2 (ii)]{su}.

\bcor\label{cellulartorickunnethtoptrivialH} We have the
following isomorphism: \be\label{topkunnethtrivpicewise}{K}^0(\mathcal{E}(X))\simeq
{K}^0(\mathcal{B})\otimes_{R(T_{comp})}PLP(\Delta).\ee Here ${K}^0(\mathcal{B})$
is a $R(T_{comp})$-module via the map which takes any
$[\mathbf{V}]\in R(T_{comp})$ for a $T$-representation $\mathbf{V}$ to
$[\mathcal{E}\times_{T} \mathbf{V}]\in {K}^0(\mathcal{B})$.  Further, we have
$K^{-1}(\mathcal{E}(X))=0$. \ecor

In the case when $X(\Delta)$ is a smooth, $PLP(\Delta)$ is isomorphic
to the Stanley-Reisner ring of $\Delta$, which we denote by
$SR(\Delta)$ (see \cite{vv}, \cite{baggio}). In
the Laurent polynomial ring
$R(T_{comp})[X_1^{\pm1},\ldots, X_d^{\pm1}]$, let
$\displaystyle X_{F}:=\prod_{\rho_i\in F}(1-X_i)$ for every
$F\subseteq \Delta(1)$.  As an $R(T_{comp})$-algebra, $SR(\Delta)$ has
the following presentation:
$R(T_{comp})[X_1^{\pm1},\ldots, X_d^{\pm1}]/\mathcal{I}$ where
$\mathcal{I}$ is the ideal generated by the elements $X_{F}$ where
$\{v_j~\mid ~\rho_j\in F\}$ do not span a cone in $\Delta$ and
$\displaystyle\left(\prod_{\rho_j\in \cf_{+}(1)}X_j^{\langle
    u,v_j\rangle}\right)-e^{u}$ for $u\in M$.

We therefore have the following corollary when $X(\Delta)$ is a smooth
complete $T$-cellular toric variety.

\bcor\label{eqkthtbsmooth} As a
${K}_{H_{comp}}^0(\mathcal{B})$-algebra the ring
${K}_{H_{comp}}^0(\mathcal{E}(X))$ has the following presentation:
$\mathcal{R}({K}_{H_{comp}}^0(\mathcal{B}), \Delta):=
{K}_{H_{comp}}^0(\mathcal{B})[X_1^{\pm1},\ldots,
X_d^{\pm1}]/\mathcal{J}$ where $\mathcal{J}$ is the ideal generated by
the elements $X_{F}$ where $\{v_j~\mid ~\rho_j\in F\}$ do not span a
cone in $\Delta$ and
$\displaystyle \left(\prod_{\rho_j\in \cf_{+}(1)}X_j^{\langle
    u,v_j\rangle}\right)-[\mathcal{L}_u]$ for $u\in M$.\ecor

\section{Equivariant $K$-theory of horospherical
  varieties}\label{horospherical}

We refer to \cite{pasquier}, \cite[Chapter 13]{perrin} and
\cite{monahan} for basic results in the theory of horospherical
varieties (also see \cite{knop}, \cite{t} and \cite{perrin1}).

Let $G$ be a connected complex reductive linear algebraic group, let
$B\subseteq G$ be a Borel subgroup, $T\subseteq B$ a maximal torus. We
denote by $\Phi$ the system of roots of $(G, T )$, by $\Phi_+$ the
system of positive roots with respect to $T\subseteq B$, and by
$\Sigma$ the system of simple roots. Let $W$ denote the Weyl group of
$(G, T )$.  For $I \subseteq \Sigma$, we denote by $W_I$ the subgroup
of $W$ generated by the simple reflections $s_{\alpha}$ for
$\alpha \in I$. We denote by $P_I$ the parabolic subgroup of $G$
generated by $B$ and $W_{I}$. There is a one-one correspondence
between the subsets $I\subseteq \Sigma$ of simple roots and the
parabolic subgroups $P_{I}$ of $G$ containing $B$.

A {\it homogeneous space} ${G}/{H}$ is called {\it horospherical} if
${H}$ contains a maximal unipotent subgroup ${U}$
of ${G}$. We shall assume that $U$ is the unipotent radical of $B$.

An {\it embedding} of a homogeneous space $G/H$ is a pair $(X, x)$,
where $X$ is a normal $G$-variety and $x$ is a point of $X$, such that
the orbit of $x$ in $X$ is open and isomorphic to $G/H$. When $G/H$ is
horospherical we call $X$ a horospherical variety.

Recall that a homogeneous space $G/H$ is called {\em spherical} if it
contains an open orbit under the action of a Borel subgroup $B$ of
$G$. A spherical variety is an embedding of a spherical homogeneous
space $G/H$. By Bruhat decomposition it follows that every
horospherical variety is spherical.

For a spherical homogeneous space $G/H$ we denote by $D$ the
collection of irreducible divisors of $G/H$ which are $B$-stable but
not $G$-stable. The elements of $D$ are called {\it colors}.  Let $X$
be an embedding of $G/H$. We denote by $X_1, \ldots , X_m$ the
$G$-stable irreducible divisors of $X$. We can identify $D$ with the
set of irreducible divisors of $X$ which are $B$-stable and not
$G$-stable. Thus $D ~\cup ~\{X_1, \ldots , X_m\}$ is the set of
irreducible $B$-stable divisors of $X$.

A {\it color of $X$} is a  color which contains a closed  $G$-orbit.

A spherical variety $X$ is {\it toroidal} if it does not contain any
color.

Let $G/H$ be a horospherical homogeneous space.  Let
${P}:=N_{{G}}({H})$. Then ${P}$ is a parabolic subgroup of ${G}$ and
${P}/{H}$ is a torus. Consider the ${P}/{H}$-bundle
${G}/{H}\ra {G}/{P}$.

Let ${Y}$ be a ${P}/{H}$-toric variety.  We can consider the
associated bundle ${X}:={G}\times_{{P}} {Y}$ where ${P}$ acts by left
multiplication on ${Y}$ through the projection ${P}/{H}$.  Then $X$
contains $G/H=G\times_{P} P/H$ as a dense orbit so that it is an
embedding of the horospherical homogeneous space $G/H$. It is further
known that $X$ is toroidal (see \cite{pasquier}). Thus
$X={G}\times_{{P}} {Y}$ is a toroidal horospherical variety.

It is in fact known that any toroidal horospherical embedding of $G/H$
is of the form $G\times_{P} Y$ with $Y$ a toric variety for the torus
$P/H$.  Indeed, $Y$ is the toric variety corresponding to the fan
$\mathbb{F}(X)$ associated to $X$ (see \cite{pasquier} and
\cite{perrin}).

Let $G_{comp}$ denote a maximal compact subgroup of $G$ such that
$B_{comp}=B\cap G_{comp}=T_{comp}$ is a maximal torus of $G_{comp}$
and $Z:=P\cap G_{comp}$ is the centralizer of a circle subgroup of
$G_{comp}$.

Moreover, $H_{comp}=H\cap G_{comp}$ is a compact subgroup of $H$ such
that $N_{G_{comp}} (H_{comp})=Z$. Thus $(P/H)_{comp}\simeq Z/H_{comp}$
and $G/P\simeq G_{comp}/Z$.

Let $R(Z)$ denote the representation ring of $Z$. Let $R(Z/H_{comp})$
denote the representation ring of $Z/H_{comp}$. The
projection $Z\ra Z/H_{comp}$ induces a canonical restriction map from
\[R(Z/H_{comp})=R(P/H)\stackrel{r}{\ra} R(Z)=R(P).\]

Recall that $K^0_{G_{comp}}(G/P)$ is an $R(Z)$-algebra via the map
which sends any $[V]\in R(Z)$ to the class of the associated
$G_{comp}$-equivariant vector bundle $[G_{comp}\times_{Z} V]$ where
$G_{comp}$-acts from the left. Thus $K^0_{G_{comp}}(G/P)$ gets the
structure of $R(Z/H_{comp})$-algebra by composing with $r$.
  
Recall that $X$ is complete if and only if $\mathbb{F}(X)$ is a
complete fan (see \cite[Theorem 4.2]{knop},\cite[Proposition
5.3.11]{monahan}). Further, recall that $X$ is smooth if and only if
the fan $\mathbb{F}(X)$ is smooth (see \cite[Theorem 4.2.3]{perrin1}
and \cite[Proposition 6.3.2]{monahan}).


Suppose that $Y$ is a complete $P/H$-cellular toric variety. Then we
have the following result:

  \bth\label{equivhorospherical} The $G_{comp}$-equivariant
  topological $K$-ring of a complete toroidal horospherical embedding
  $X$ of $G/H$ is isomorphic to the following ring
  $K^0_{G_{comp}}(G/P)\otimes_{R((P/H)_{comp})}PLP(\mathbb{F}(X))$ as
  an $K^0_{G_{comp}}(G/P)$-algebra.

  \eeth \begin{proof} Recall that we have a $G$-equivariant
    isomorphism of $X$ with $G\times_{P} Y$. Further, note that
    $G\times_{P} Y=\Big(G\times_{P} P/H\Big)\times_{P/H} Y$. Thus we
    have a $G$-equivariant isomorphism of $X$ with the $G$-equivariant
    toric bundle $\mathcal{E}\times_{P/H}Y$ associated to the
    $G$-equivariant principal $P/H$-bundle:
    $\mathcal{E}=G/H\ra \mathcal{B}=G/P$.

    Recall that $G/P$ is a smooth projective variety and has a
    filtrable Bialynicki Birula $T$-cellular decomposition given by
    the Schubert cells indexed by the minimal length coset
    representatives $W^{I}$ of $W/W_{I}$, where $P=P_I$ for
    $I\subseteq \Sigma$ (see \cite{kk} and \cite[Proposition
    5.1]{BGG}). Hence by Proposition \ref{celldec}, it follows that
    $K_{T_{comp}}^0(G/P)$ is a free $R(T_{comp})$-module of rank
    $|W/W_{I}|$ and $K_{T_{comp}}^{-1}(G/P)=0$. By \cite[Theorem
    4.4]{ml}, this further implies that
    $K^{-1}_{G_{comp}}(G/P)=K^{-1}_{T_{comp}}(G/P)^{W}=0$. Thus
    Assumption \ref{base} holds for $\mathcal{B}=G/P$. Since $Y$ is a
    complete $P/H$-cellular toric variety, by Theorem \ref{main2},
    $K^0_{(P/H)_{comp}}(Y)$ is isomorphic as an
    $R((P/H)_{comp})$-algebra to $PLP(\mathbb{F}(X))$. The theorem now
    follows from Theorem \ref{cellulartorickunnethtop}, by taking
    $\mathcal{E}=G/H$, $\mathcal{B}=G/P$, $T=P/H$ and
    $\Delta=\mathbb{F}(X)$.
  \end{proof}

  We have the following description of $K^0_{G_{comp}}(X)$ in the case
  when $Y$ is a smooth complete $P/H$-cellular toric variety.  See
  \cite[Section 3.2]{kh} for the analogous result on the ordinary
  cohomology ring of a smooth complete horospherical variety.

  \bcor\label{smoothcase} The $G_{comp}$-equivariant topological
  $K$-ring of a smooth complete toroidal horospherical embedding $X$
  of $G/H$, is isomorphic to the ring
  $K^0_{G_{comp}}(G/P)\otimes_{R((P/H)_{comp})}SR(\mathbb{F}(X))$ as an
  $K^0_{G_{comp}}(G/P)$-algebra. \ecor

  \begin{proof}
    Since $Y$ is a smooth $P/H$-cellular toric variety,
    $K_{(P/H)_{comp}}^0(Y)$ is isomorphic to $\mathcal{K}^0_{P/H}(Y)$
    which in turn is isomorphic to $SR(\mathbb{F}(X))$, namely the
    Stanley-Reisner ring of the fan $\mathbb{F}(X)$, as an
    $R(P/H)$-algebra (see \cite{vv}, \cite{baggio}).  The proof now
    follows from Theorem \ref{equivhorospherical} and Corollary
    \ref{eqkthtbsmooth}.
  \end{proof}

  
  We get the following description of the $G$-equivariant Grothendieck
  ring of $X$ when $Y$ is a smooth complete $P/H$-cellular toric
  variety. The corresponding result for the ordinary Grothendieck ring
  can alternately be obtained by applying \cite[Proposition 3.1]{u3}.

  \bcor\label{smoothcaseGroth} (1) The $G$-equivariant Grothendieck ring
  $\mathcal{K}^0_{G}(X)$ of a smooth complete toroidal horospherical
  embedding $X$ of $G/H$ is isomorphic to
  \[\mathcal{K}^0_{G}(G/P)\otimes_{R(P/H)}SR(\mathbb{F}(X))\] as an
  $\mathcal{K}^0_{G}(G/P)$-algebra.  

  (2) We have the presentation
  $\mathcal{K}^0_{G}(X)\cong \mathcal{R}(
  \mathcal{K}^0_{G}(G/P),\mathbb{F}(X)))$ as an
  $\mathcal{K}^0_{G}(G/P)$-algebra.
  \ecor
  \begin{proof}
    Since $G/P$ is smooth and projective, we have \be\label{iso1}
    \mathcal{K}^0_{G}(G/P) \cong K^0_{G_{comp}}(G/P).\ee Thus when $Y$
    is smooth complete and $P/H$-cellular we have \be\label{iso4}
    \mathcal{K}^0_{G}(\mathcal{E}\times_{P/H} Y) \cong
    K^0_{G_{comp}}(\mathcal{E}\times_{P/H} Y) \ee (see
    \cite[Proposition 5.5.6]{cg}). Again since $Y$ is a smooth
    complete $P/H$-cellular toric variety we have
    \be\label{iso2}\mathcal{K}^0_{P/H}(Y)\cong
    K^0_{(P/H)_{comp}}(Y)\cong SR(\mathbb{F}(X)).\ee Now,
    \be\label{iso3} R(P/H)= R(Z/H_{comp}).\ee Recall that
    $\mathcal{K}^0_{G}(G/P)$ is an $R(P)$-algebra via the map which
    sends any $[V]\in R(P)$ to the class of the associated
    $G$-equivariant vector bundle $[G\times_{P} V]$ where $G$-acts
    from the left. Thus $\mathcal{K}^0_{G}(G/P)$ gets the structure of
    an $R(P/H)$-algebra by composing with $r:R(P/H)\ra R(P)$.
    Furthermore, \eqref{iso1} and \eqref{iso2} are isomorphisms of
    $R(P/H)=R(Z/H_{comp})$-algebras and \eqref{iso4} is an isomorphism
    of $R(G)=R(G_{comp})$-algebras. Moreover, since there is a
    $G$-equivariant isomorphism of $X$ with the $G$-equivariant toric
    bundle $\mathcal{E}\times_{P/H} Y$ we have \be\label{iso5}
    \mathcal{K}^0_{G}(X)
    \cong\mathcal{K}^0_{G}(\mathcal{E}\times_{P/H} Y) .\ee The proof
    follows successively from \eqref{iso5}, \eqref{iso4}, Corollary
    \ref{eqkthtbsmooth}, \eqref{iso1}, \eqref{iso2} and \eqref{iso3}.



\end{proof}

\noindent {\bf Acknowledgments:} I am very grateful to Prof. Michel
Brion for several valuable comments and suggestions on earlier
versions of this manuscript. I thank Prof. Parameswaran Sankaran for
going through the manuscript and for valuable comments.


\noindent {\bf Funding:} I thank SERB MATRICS grant no:
MTR/2022/000484 for financial support.

\noindent{\bf Conflict of Interest statement:} The author has no
relevant financial or non-financial interests to disclose.

\end{document}